\documentclass[12pt]{article}
\usepackage{amsmath,amsthm,amsfonts,amscd,amssymb} 
\usepackage[T1]{fontenc}
\usepackage[utf8]{inputenc}
\usepackage{authblk}
\usepackage{verbatim}   
\usepackage{eucal} 
\usepackage{longtable} 
\usepackage{makeidx}
\usepackage{url}
\usepackage[all]{xy} 
\usepackage{enumerate}

\newcommand{\field}[1]{\mathbb{#1}}

\newcommand{\F}{\field{F}}

\newfont{\cyrr}{wncyr10}

\newcommand{\thml}[1]{\begin{theorem}\label{#1}}
\newcommand{\mht}{\end{theorem}}
\newcommand{\cnj}{\begin{conjecture}}
\newcommand{\cnjl}[1]{\begin{conjecture}\label{#1}}
\newcommand{\jnc}{\end{conjecture}}
\newcommand{\dfnl}[1]{\begin{definition}\label{#1}}
\newcommand{\nfd}{\end{definition}}
\newcommand{\pro}{\begin{proposition}}
\newcommand{\prol}[1]{\begin{proposition}\label{#1}}
\newcommand{\orp}{\end{proposition}}
\newcommand{\crl}{\begin{corollary}}
\newcommand{\crll}[1]{\begin{corollary}\label{#1}}
\newcommand{\lrc}{\end{corollary}}
\newcommand{\lmm}{\begin{lemma}}
\newcommand{\lmml}[1]{\begin{lemma}\label{#1}}
\newcommand{\mml}{\end{lemma}}
\newcommand{\prf}{\begin{proof}}
\newcommand{\prfl}[1]{\begin{proof}\label{#1}}
\newcommand{\frp}{\end{proof}}
\newcommand{\axi}{\begin{axiom}}
\newcommand{\axil}[1]{\begin{axiom}\label{#1}}
\newcommand{\ixa}{\end{axiom}}
\newcommand{\rmkl}[1]{\begin{remark}\label{#1}}
\newcommand{\kmr}{\end{remark}}
\newcommand{\exa}{\begin{example}}
\newcommand{\exal}[1]{\begin{example}\label{#1}}
\newcommand{\axe}{\end{example}}
\newcommand{\alg}{\begin{algorithm}}
\newcommand{\gla}{\end{algorithm}}
\newcommand{\nte}{\begin{note}}
\newcommand{\ntel}[1]{\begin{note}\label{#1}}
\newcommand{\etn}{\end{note}}
\newcommand{\app}{\begin{application}}
\newcommand{\appl}[1]{\begin{application}\label{#1}}
\newcommand{\ppa}{\end{application}}
\newcommand{\mat}{\begin{matrix}}
\newcommand{\tam}{\end{matrix}}
\newcommand{\smm}{\begin{summary}}
\newcommand{\mms}{\end{summary}}

\newcommand{\teo}{\begin{teorema}}
\newcommand{\teol}[1]{\begin{teorema}\label{#1}}
\newcommand{\oet}{\end{teorema}}
\newcommand{\df}{\begin{definicion}}
\newcommand{\dfl}[1]{\begin{definicion}\label{#1}}
\newcommand{\fd}{\end{definicion}}
\newcommand{\por}{\begin{proposicion}}
\newcommand{\porl}[1]{\begin{proposicion}\label{#1}}
\newcommand{\rop}{\end{proposicion}}
\newcommand{\corl}[1]{\begin{corolario}\label{#1}}
\newcommand{\roc}{\end{corolario}}
\newcommand{\leml}[1]{\begin{lema}\label{#1}}
\newcommand{\mel}{\end{lema}}
\newcommand{\pru}{\begin{prueba}}
\newcommand{\prul}[1]{\begin{prueba}\label{#1}}
\newcommand{\urp}{\end{prueba}}
\newcommand{\axa}{\begin{axioma}}
\newcommand{\axal}[1]{\begin{axioma}\label{#1}}
\newcommand{\xa}{\end{axioma}}
\newcommand{\nta}{\begin{nota}}
\newcommand{\ntal}[1]{\begin{nota}\label{#1}}
\newcommand{\atn}{\end{nota}}
\newcommand{\eje}{\begin{ejemplo}}
\newcommand{\ejel}[1]{\begin{ejemplo}\label{#1}}
\newcommand{\je}{\end{ejemplo}}
\newcommand{\api}{\begin{aplicacion}}
\newcommand{\apil}[1]{\begin{aplicacion}\label{#1}}
\newcommand{\ipa}{\end{aplicacion}}

\newcommand{\ali}{\begin{align}}
\newcommand{\ila}{\end{align}}
\newcommand{\enu}{\begin{enumerate}}
\newcommand{\une}{\end{enumerate}}
\newcommand{\arr}{\begin{array}}
\newcommand{\rra}{\end{array}}
\newcommand{\eqa}{\begin{eqnarray}}
\newcommand{\aqe}{\end{eqnarray}}
\newcommand{\equ}{\begin{equation}}
\newcommand{\uqe}{\end{equation}}
\newcommand{\subq}{\begin{subequations}}
\newcommand{\qbus}{\end{subequations}}

\newtheorem{theorem}{Theorem}[section]
\newtheorem{lemma}[theorem]{Lemma}
\newtheorem{proposition}[theorem]{Proposition}
\newtheorem{corollary}[theorem]{Corollary}
\newtheorem{conjecture}[theorem]{Conjecture}

\theoremstyle{definition}
\newtheorem{definition}{Definition}[section]
\newtheorem{axiom}[definition]{Axiom}

\theoremstyle{remark}
\newtheorem{remark}{Remark}[section]
\newtheorem{example}[remark]{Example}
\newtheorem{application}[remark]{Application}
\newtheorem{algorithm}[remark]{Algorithm}
\newtheorem{note}[remark]{Note}
\newtheorem{summary}[remark]{Summary}

\newtheorem{teorema}[theorem]{Teorema}
\newtheorem{lema}[theorem]{Lema}
\newtheorem{proposicion}[theorem]{Proposici\'on}
\newtheorem{corolario}[theorem]{Corolario}

\newtheorem{definicion}[definition]{Definici\'on}
\newtheorem{axioma}[definition]{Axioma}

\newtheorem{nota}[remark]{Nota}
\newtheorem{ejemplo}[remark]{Ejemplo}
\newtheorem{aplicacion}[remark]{Aplicaci\'on}

\author[1]{
Victor Bautista-Ancona}
\author[1]{Javier Diaz-Vargas}
\author[1]{Jos\'e Alejandro Lara Rodr\'{\i}guez}
\author[2]{Francisco X. Portillo-Bobadilla}
\affil[1]{Universidad Aut\'onoma de Yucat\'an}
\affil[2]{Universidad Aut\'onoma de la Ciudad de M\'exico}

\title{A description of a Drinfeld module with class number $h=1$ and rank $1$}

\begin{document}
\maketitle

\abstract{We work with detail the Drinfeld module over the ring
$$A=\F_2[x,y]/(y^2+y=x^3+x+1).$$  The example in question is one of the
four examples that come from quadratic imaginary fields with class number $h = 1$
and rank one.

We develop specific formulas for the coefficients $d_k$ and $\ell_k$ of the
exponential and logarithmic functions and relate them with the product $D_k$ of
all monic elements of $A$ of degree $k$. On the Carlitz module, $D_k$ and $d_k$
coincide, but this is not true in general Drinfeld modules. On this example, we
obtain a formula relating both invariants. We prove also using elementary
methods a theorem due to Thakur that relate two different combinatorial symbols
important in the analysis of solitons.

}

\section{Introduction.}

Let $\mathbb{F}_q$ be a finite field of characteristic $p$ and $K$ a
function field over $\mathbb{F}_q$. After we choose $\infty$, a
fixed infinite place of $K$, let $A$ be the ring of regular
functions outside of $\infty$ and let $K_{\infty}$ be its
completion. Now take $\mathbb{C}_{\infty}$ to be the completion of
an algebraic closure of $K_{\infty}$.

Let $\mathbb{C}_\infty\{\tau\}$ be the ring of \emph{twisted polynomials},
i.e., the noncommutative ring of polynomials $\sum a_i \tau^i$ with
coefficients in $\mathbb{C}_\infty$ such that $\tau z = z^q \tau$.
A twisted polynomial $f = a_0 + a_1 \tau + \dotsb + a_d \tau^d\in
\mathbb{C}_\infty$ is identified with the $\mathbb{F}_q$-linear endomorphism of
$\mathbb{C}_\infty$,
\begin{align*}
z \mapsto f(z) = a_0 z+ a_1z^q + \dotsb + a_d z^{q^d}.
\end{align*}

A Drinfeld $A$-module of rank one is a $\mathbb{F}_q$-algebra homomorphism
$\rho \colon  A \to \mathbb{C}_{\infty}\{ \tau_p\}$ injective, for which
$\rho(a) = a \tau^0 + \text{higher order terms in } \tau$.
The action $a \cdot z =
\rho(a)(z)$ of $A$ in $\mathbb{C}_\infty$ makes $\mathbb{C}_\infty$ into an
$A$-module, and hence the name ``Drinfeld module''.

For each Drinfeld module $\rho$ we associate an exponential  entire function $e$ defined for a
power series in all $\mathbb{C}_{\infty}$ by
    \[
        e(z) = \sum_{i=0}^{\infty} \frac{z^{q^i}}{d_i}.
    \]
The linear term in this exponential function satisfy the following fundamental
functional equation
    \begin{equation}
        e(az) = \rho_a (e(z)), \label{eq:functionaleq}
    \end {equation}
for $z \in \mathbb{C}_{\infty}$ and $a \in A$, where $\rho_a$ stands
for $\rho(a)$.

The  Carlitz module, defined by Carlitz \cite{C1} in 1935, is given
by the $\mathbb{F}_q$-algebra homomorphism $C : \mathbb{F}_q [t] \to
\mathbb{C}_\infty\{\tau\}$ determined by $C_t = t+\tau^q$. Equation
(\ref{eq:functionaleq}) produces $e(tz)=t e(z)+e(z)^q$. It follows
that
        \[
            \sum_{i=0}^{\infty} \frac{(t^{q^i}-t)z^{q^i}}{d_i}
= \sum_{i=0}^{\infty} \frac{ z^{q^{i+1}}}{d_i^q}
        \]
By equating coefficients we get a unique solution $d_n = [n] d_{n-1}^q$ where
$[n] = (t^{q^n} - t)$ and $d_0 = 1$. Therefore, $d_n = [n][n-1]^q \dotsm
[1]^{q^{n-1} }$ and it is easily seen that $d_n$ is the product of all monic
polynomials of degree $n$.

Since $e(z)$ is periodic, it can not have a global inverse, but we
may formally derive an inverse $\log (z)$ for $e(z)$ as a power
series around the origin. By definition $e(\log (z))=z$. Since
$e(z)$ satisfies the functional equation  $e(tz) = te(z)+ e(z)^q$,
it follows that $tz = \log (te(z)) + \log (e(z)^q)$. Replacing $\log
(z)$ for $z$ we obtain $t \log (z) = \log (tz) + \log (z^q)$. Let
$\log (z) = \sum z^{q^i}/\ell_i$. Then
\begin{align*}
  \sum _{i=0}^\infty \frac{(t-t^{q^i})\cdot z^{q^i} }{\ell_i}=
\sum _{i=0}^\infty \frac{z^{q^{i+1} } }{\ell_i}
\end{align*}
It follows that $\ell_{i+1} = -[i+1] \ell_i$. Therefore $\ell_i = (-1)^i
[i][i-1] \dotsm [1]$.

We follow the ideas developed in the Carlitz module case, but applied to the Drinfeld module over
$A=\F_2[x,y]/(y^2+y=x^3+x+1)$. We explore specific ways to understand the mentioned example,
which is one of four examples provided from imaginary quadratic fields with class number h = 1
\cite{mcrae} and rank 1. The formulas obtained are compared with the Theorem $4.15.4$ of \cite{thakur}
and are related to solitons, as exposed in Chapter $8$ of the same reference, and Theorem 3 of the article \cite{thakur2}.

\section{Action of the Drinfeld module on the variables $x$ and $y$.}

In our example, we have $d_\infty =1 $, $v_\infty (x)=-2$, $v_\infty (y)=-3$,
and using that $\deg(a)=-v_\infty(a) d_\infty$ $\forall a\in A$, it follows that
$\deg(x)=2$ and $\deg(y)=3$.

Based on it, the Drinfeld Module $\rho$ is determined by its values in $x$ and $y$
(actually, it is enough to know its value in one element $a\in A$, see 2.5 in
\cite{thakur}).
 According to the aforementioned degrees and that the unique sign in our example
is $+1$, we obtained that
\begin{align*}
\rho_x & = x+x_1\tau+\tau^2,\\
\rho_y &=y+y_1\tau+y_2\tau^2+\tau^3
\end{align*}
with $x_1$, $y_1$, $y_2\in A$.
Now, using the commutative property of the Drinfield module $\rho_x\rho_y=\rho_y\rho_x$ and equaling on degree $1$,
we get $$x_1(y^2+y)=y_1(x^2+x).$$
Next, using the equation on the curve $y^2+y=x^3+x+1$ and dividing, we obtain
$$y_1=x_1\left(x+1+\frac{1}{x^2+x}\right).$$
This implies that  $x^2+x\mid x_1$ and $y^2+y\mid y_1$.
Assuming that $x_1=x^2+x$, it is also obtained that $y_1=y^2+y$.
Now, equaling on degree $2$, one has the equation
\equ\label{grado2}
(x^4+x)y_2=-y_1x_1^2+y_1^2x_1+(y^4+y).
\uqe
But, we can use the identities
\begin{align*}
y^4+y &=(y^2+y)^2+y^2+y\\
&=(y^2+y)(y^2+y+1)\\
&=(y^2+y)(x^3+x)\\
&=(y^2+y)(x^2+x)(x+1)
\end{align*}
%
and
$$x^4+x=(x^2+x)(x^2+x+1).$$
So dividing the equation (\ref{grado2}) by $x_1=x^2+x$, and substituting the
values $x_1$ and $ y_1 $, we  get
\begin{align*}
y_2(x^2+x+1)&=(y^2+y)(x^2+x+y^2+y)+(y^2+y)(x+1)\\
&=(y^2+y)(y^2+y+x^2+1)\\
&=(y^2+y)(x^3+x^2+x).
\end{align*}
%
Thus, clearing $y_2$,  we have $y_2=x(y^2+y)$, as it is known in the literature
\cite{hayes}.

\section{Exponential and Logarithm coefficients.}

We find recursive formulas for the coefficients of both the exponential $e(z)$ and the logarithmic $\log(z)$ functions asociated to Drinfeld module in $ A $.

Write $$e(z)=\sum_{i=0}^\infty \frac{z^{2^i}}{d_i}=\sum_{i=0}^\infty
a_iz^{2^i}$$
and
$$\log(z)=\sum_{i=0}^\infty \frac{z^{2^i}}{\ell_i}=\sum_{i=0}^\infty
b_iz^{2^i}$$
where $a_i=d_i^{-1}$ y $b_i=\ell_i^{-1}$.
Using that
\begin{align*}
e(xz) &=\rho_x\left(e(z)\right)\\
& =xe(z)+[1]_xe^2(z)+e^4(z)
\end{align*}
%
where $[1]_x=x^2+x$.
Then, working both sides of the equality:
$$e(xz)+xe(z)=[1]_xe^2(z)+e^4(z),$$
we have on the left side:
\equ\label{lado_izquierdo}
\begin{split}
e(xz)+xe(z) &= \sum_{j=0}^\infty \left(x^{2^j}+x\right)a_j z^{2^j}\\
&= \sum_{j=0}^\infty [j]_x a_j z^{2^j}\\
&= [1]_x a_1 z^2+\sum_{j=2}^\infty [j]_x a_j z^{2^j},
\end{split}
\uqe
where $[j]_x:=x^{2^j}+x$.
Now, developing the right side, we get:
\begin{align*}
[1]_xe^2(z)+e^4(z)=[1]_x\sum_{i=0}^\infty a_i^2
z^{2^{i+1}}+\sum_{i=0}^\infty a_i^4 z^{2^{i+2}}.
\end{align*}
%
From where, by setting $j=i+1$ in the first sum, and $j=i+2$ in the second sum, we obtain:
\equ
\begin{split}\label{lado_derecho}
[1]_xe^2(z)+e^4(z) &=[1]_x\sum_{j=1}^\infty a_{j-1}^2
z^{2^{j}}+\sum_{j=2}^\infty a_{j-2}^4 z^{2^{j}}\\
&= [1]_x a_0^2 z^2+\sum_{j=2}^\infty \left([1]_x a_{j-1}^2+a_{j-2}^4\right)
z^{2^{j}}.
\end{split}
\uqe
Comparing equations (\ref{lado_izquierdo}) and (\ref{lado_derecho}), recursive
formulas are obtained
\begin{align}
a_1 &= a_0^2 \nonumber\\
a_j &=\frac{[1]_x a_{j-1}^2+a_{j-2}^4}{[j]_x} \quad \mbox{for $j\geq
2$.}\label{recursive_aes}
\end{align}
%
Subsequently, we assume that $a_0=1$, i.e., the exponential is normalized.
Notice that if we do not normalize the coefficients, the exponential function varies by a factor given by the initial term. If we denote $e(z,a_0)$ to this exponential function, it is easy to see that
\equ
e(z,a_0)=a_0 e(z),
\uqe
where $e(z)$ is the normalized exponential.

Now, in terms of the  $ d_j$'s  (assuming also, the normalization of the
exponential), the recursive formula is as follows:
\begin{align}
d_1 &=d_0^2=1  \nonumber \\
d_j &=\frac{[j]_xd_{j-1}^2d_{j-2}^4}{[1]_xd_{j-2}^4+d_{j-1}^2} \quad\mbox{for
$j\geq
2$.}\label{recursive_Dees}
\end{align}
%
Similarly, for the logarithm function, we have that
\begin{align*}
x\log(z) &= \log\left(\rho_x(z)\right)\\
&=\log\left(xz+[1]_xz^2+z^4\right)\\
&=\log(xz)+\log([1]_xz^2)+\log(z^4),
\end{align*}
%
from which it follows that $$x\log(z)+\log(xz)=\log([1]_x z^2)+\log(z^4).$$
So, we developed the left side to
\equ\label{log_left}
x\log(z)+\log(xz)=\sum_{j=0}^\infty(x^{2^j}+x)b_j
z^{2^j}=\sum_{j=1}^\infty[j]_x b_j z^{2^j}.
\uqe
Note that $[0]_x=0$. The right side must be
$$\log([1]_x z^2)+\log(z^4)=\sum_{i=0}^\infty [1]_x^{2^i} b_i z^{2^{i+1}}+\sum_{i=0}^\infty b_i z^{2^{i+2}}.$$
Again, by setting $ j = i + 1 $ in the first sum, and $ j = i + 2 $ in the second sum, we obtain
\equ\label{log_right}
\log([1]_x z^2)+\log(z^4)=[1]_x b_1 z^2+\sum_{j=2}^\infty
\left([1]_x^{2^{j-1}}b_{j-1}+b_{j-2}\right)z^{2^j}.
\uqe
%
Comparing the terms in the equations (\ref{log_left}) and (\ref{log_right}),
we obtain the recursive formulas:
\begin{align}
b_1 & = b_0 \nonumber\\
b_j &= \frac{[1]_x^{2^{j-1}}b_{j-1}+b_{j-2}}{[j]_x} \quad\mbox{for $j\geq
2$.}\label{recursive_bees}
\end{align}
Now again, if $\log(z,b_0)$ is the logarithmic function with initial term $b_0$, and $\log(z)=\log(z,1)$ is the normalized logarithm, by the recursion formula, we deduce the relation:
\equ
    \log(z,b_0)=b_0\log(z).
\uqe

In terms of values $\ell_i$'s, the recursions are as follows:
\begin{align*}
\ell_1 &=\ell_0 \nonumber\\
\ell_j &= \frac{[j]_x \ell_{j-1}\ell_{j-2}}{[1]_x^{2^{j-1}}\ell_{j-2}+
\ell_{j-1}}
\quad\mbox{for $j\geq
2$}.
\end{align*}
%

\section{Formul\ae\mbox{ }for computing $\rho_a$.}

The first formula is recursive and is in the spirit of the proposition $ 3.3.10 $ in \cite{goss}.

Assume that $\rho_a=\sum_{k=0}^d \rho_{a,k} \tau^k$ with $d=\deg(a)$.
We will use again commutativity $\rho_x\rho_a=\rho_a\rho_x$ and the explicit expression: $\rho_x=x+[1]_x\tau+\tau^2$.
Then, multiplying
\begin{align*}
\rho_x\rho_a &=(x+[1]_x\tau+\tau^2)\left(\sum_{k=0}^d \rho_{a,k} \tau^k\right)\\
&=\sum_{k=0}^d
\left(x\rho_{a,k}\tau^k+[1]_x\rho_{a,k}^2\tau^{k+1}+\rho_{a,k}^4\tau^{k+2}
\right)
\end{align*}
%
and multiplying
\begin{align*}
\rho_a\rho_x &=\left(\sum_{k=0}^d \rho_{a,k} \tau^k\right)(x+[1]_x\tau+\tau^2)\\
&=\sum_{k=0}^d\left(x^{2^k}\rho_{a,k}\tau^k+[1]_x^{2^k}\rho_{a,k}\tau^{k+1}
+\rho_{a,k}\tau^{k+2} \right).
\end{align*}
By comparing terms a recursive formula is obtained
\begin{align*}
\rho_{a,0}& =a  & & \text{(first term in recursion)} \\
\rho_{a,1} &= a^2+a & &  \text{(comparing degree $k = 1$)}
\end{align*}
$$
\rho_{a,k}=\frac{[1]_x^{2^{k-1}}\rho_{a,k-1}+\rho_{a,k-2}}{[k]_x}+\frac{[1]
_x\rho^2_{a,k-1}+\rho^4_{a,k-2}}{[k]_x},\mbox{ for $k\geq 2$}.
$$
%

Note the similarity to the recursive formulas for  $a_j$'s  and $b_j$'s  in
the previous section, equations (\ref{recursive_aes}) and
(\ref{recursive_bees}). The same phenomenon occurs in the Carlitz module, but in
such a case, there is only a single summand.



Another way to calculate $\rho_a$, is based on the use of the exponential and the logarithm functions and their formal development as power series. We know that
$$
e(a \log(z))=\rho_a(e(\log(z))=\rho_a(z).
$$
Using power series as in the previous section, we get to
\begin{align*}
\rho_a(z)&=\sum_{k=0}^\infty \left( \sum_{j=0}^k a_jb_{k-j}^{2^j}a^{2^j}\right)z^{2^k}\\
&=\sum_{k=0}^\infty \left( \sum_{j=0}^k \frac{a^{2^j}}{d_j\ell_{k-j}^{2^j}}
\right)z^{2^k}.
\end{align*}
The combinatorial terms in the sum, are the ones that D. Thakur used to develop his alternative perspective on solitons \cite{thakur2}.

We introduce the following notation for the following
pages:$$p_k(w):=\left\{\begin{matrix} w \cr
q^k\end{matrix}\right\}:=\sum_{j=0}^k \frac{w^{2^j}}{d_j \ell_{k-j}^{2^j}}. $$
Hence, using that $\rho_a=\sum \rho_{a,k} \tau^k$ is a monic polynomial in $\tau$ of
degree $\deg(a)$, we have that
\begin{align*}
p_k(a)=\rho_{a,k}=\begin{cases} 0, & \mbox{ if $\deg(a)<k$} \cr
1, & \mbox{ if $\deg(a)=k$. }\end{cases}
\end{align*}

\section{Comparing the polynomials $p_k(w)$ and $e_k(w)$.}

We define the following sets
\begin{align*}
A_{<k} &:=\{a\in A \colon \deg(a)<k\}\\
A_{k} & :=\{a\in A \colon \deg(a)=k\}
\intertext{and the polynomial}
e_k(w)&=\prod_{a\in A_{<k}} (w+a).
\end{align*}
%
Clearly, as every element of $A_{<k}$ is a root of $ p_k (w) $, we have
that
$$R_k(w):=\frac{p_k(w)}{e_k(w)}$$
is a polynomial.
In addition, as $p'_k(w)=a_k=\ell _k^{-1}\neq 0$, $p_k(w)$ and
$R_k(w)$ have
no double roots. 

In order to calculate the polynomial $R_k(w)$, suppose
\begin{align*}
p_k(w)=\sum_{i=0}^k A_{k,i}w^{2^i}
\end{align*}
and
\equ\label{sum_e_k}
e_k(w)=\sum_{i=0}^{k-1} B_{k,i}w^{2^i}.
\uqe
Then, we have the following result:

\begin{theorem}\label{division_theorem}
$R_k(w)=\frac{1}{d_k}e_k(w)+C$, where
$C=\frac{1}{d_{k-1}}+\frac{B^2_{k,k-2}}{d_k}$.
\end{theorem}

\prf

Only for the purpose of this proof, suppose $k$ is fixed and
write $A_i=A_{k,i}$ and $B_i=B_{k,i}$.
Now, directly dividing $p_k$ between $e_k$, using that $ e_k $ is monic, the first term of the quotient
ratio is $A_{k}w^{2^{k-1}}$. Then, in the first line of the waste division,
we have:
\begin{multline*}
A_{k}B_{k-2}w^{2^{k-1}+2^{k-2}}+A_{k}B_{k-3}w^{2^{k-1}+2^{k-3}}+\dotsb +\\ A_{k}
B_{0}w^{2^{k-1}+1}+A_{k-1}w^{2^{k-1}}+\mbox{ lower terms.}
\end{multline*}
This implies that the next term of the quotient is
$A_{k}B_{k-2}w^{2^{k-2}}$,
  and therefore, multiplying by the summands
of $ e_k $, after cancelation of the term $A_{k}B_{k-2}w^{2^{k-1}+2^{k-2}}$,
new summands will be incorporated into the residue in the positions corresponding to the powers:
$$w^{2^{k-1}},w^{2^{k-2}+2^{k-3}},\dotsc,w^{2^{k-2}+1}.$$
Hence, all the new
terms fall into the ``lower terms" of the waste division with exception of the coefficient on $w^{2^{k-1}}$. This coefficient is $A_{k-1}+A_{k}B^2_{k-2}$.

When continuing the division and cancelling the terms of the form $A_{k}B_{j}w^{2^{k-1}+2^{j}}$ for $j<k-2$, the terms equal or higher to $w^{2^{k-1}}$ are not affected.
This ensures that the obtained quotient is:
$$A_{k}w^{2^{k-1}}+A_{k}B_{k-2}w^{2^{k-2}}+A_{k}B_{k-3}w^{2^{k-3}}+\dotsb+A_{k}
B_ { 0 } w+A_{k-1}+A_{k}B^2_{k-2}.$$
The result follows, using that $A_{k}=d_k^{-1}$ and $A_{k-1}=d_{k-1}^{-1}$.
\frp

\section{Coefficient Formulas for $e_k(w)$.}

For $k\geq 2$, set
\begin{align*}
t_k =
\begin{cases}
x^\frac{k}{2}, & \text{if $k$ is even} \cr
yx^\frac{k-3}{2}, & \text{if $k$ is odd.}\end{cases}
\end{align*}
Now, it is clear that $\deg(t_k)=k$ and that the set $\{1,t_2,\ldots,t_{k-1}\}$
is a basis of the vector space  $A_{<k}$.
Define
$D_k:=e_k(t_k)=\prod_{a\in A_k} a$.
Thus for $k \geq 3$,
\equ\label{recursive_e_k}
\begin{split}
e_k(w) &=\prod_{a\in A_{<k}} (w+a)= \prod_{a\in A_{<k-1}} (w+a)\prod_{a\in A_{k-1}} (w+a)\cr
&= \prod_{a\in A_{<k-1}} (w+a)\prod_{a\in A_{<k-1}} (w+t_{k-1}+a)\cr
&=e_{k-1}(w)e_{k-1}(w+t_{k-1})=e^2_{k-1}(w)+D_{k-1}\cdot e_{k-1}(w).
\end{split}
\uqe
Developping the right side of the equation (\ref{recursive_e_k}), we find recursive  formulas for the coefficients  $B_{k,i}$ in (\ref {sum_e_k}):

\begin{align*}
\begin{split}
e^2_{k-1}(w)+D_{k-1}\cdot e_{k-1}(w) &= \left(\sum_{i=0}^{k-2} B_{k-1,i}w^{2^i}\right)^2+D_{k-1}\left(\sum_{i=0}^{k-2} B_{k-1,i}w^{2^i}\right)\cr
&=\sum_{i=1}^{k-1} B^2_{k-1,i-1}w^{2^i}+\sum_{i=0}^{k-2}
D_{k-1}B_{k-1,i}w^{2^i}.
\end{split}
\end{align*}
%
Indeed, we have
\begin{align*}
B_{k,0} & = D_{k-1}B_{k-1,0}=D_{k-1}D_{k-2} \cdots D_2\\
B_{k,i} &= D_{k-1}B_{k-1,i}+B^2_{k-1,i-1}\\
B_{k,k-1} & = B_{k-1,k-2}=\cdots=B_{2,1}=1.
\end{align*}
Before developing explicit formulas for the coefficients $B_{k,i}$, we introduce the following symbols:
\begin{align*}
[1]_w&= w^2+w\\
[k]_w &= w^{2^k}+w.
\end{align*}
%
It is not difficult to prove that these symbols satisfy the following:

\lmm\label{k_symbol_properties} Properties of the symbol $[k]_w$.
\begin{enumerate}[1)]
\item $[k]^{2^j}_w=[k]_{w^{2^j}}$
\item $[1]_{[k]_w}=[k]_{[1]_w}$
\item $[k]_{w_1+w_2}=[k]_{w_1}+[k]_{w_2}$
\item $[k+1]_w=[k]^2_w+[1]_w$
\item $[k]_w=\sum_{i=0}^{k-1}[1]^{2^i}_w$.
\end{enumerate}
\mml


Notice that $e_k(w)$ is a polynomial on  $[1]_w$ of degree $ 2^{k-2} $.
Set
\begin{align*}
e_k(w)=\sum_{i=0}^{k-2}T_{k,i}[1]_w^{2^i}.
\end{align*}
%
Next, we will find specific formulas for the coefficients $T_{k,i}$'s.
First, define the following functions:
\begin{align*}
S_{n,r}(x_1,x_2,\cdots,x_n)=\sum_{n\geq i_1>i_2>\cdots>i_r\geq 1}\prod_{j=1}^r
x_{i_j}^{2^{n-j+1-i_j}}.
\end{align*}

We have the following lemma:

\lmm\label{square_symetric}
Properties of the sums $S_{n,r}(x_1,x_2,\cdots,x_n)$.
\begin{enumerate}[1)]
\item $S_{n,0}(x_1,\ldots,x_{n})=1$
\item $S_{n,1}(x_1,\ldots,x_{n})=x_n+x_{n-1}^{2}+\cdots+x_1^{2^{n-1}}$
\item $S_{n+1,r}(x_1,\ldots,x_{n+1})=S^2_{n,r}(x_1,\ldots,x_n)+x_{n+1}S_{n,r-1}(x_1,\ldots,x_n)$
\end{enumerate}

\mml
\prf

The first two assertions are immediate.

For the third,
note that:
\equ\label{proof_square_symetric_1}
\begin{split}
S^2_{n,r}(x_1,\ldots,x_n)&=\left(\sum_{n\geq i_1>i_2>\cdots>i_r\geq 1}\prod_{j=1}^r x_{i_j}^{2^{n-j+1-i_j}}\right)^2\cr
&=\sum_{n\geq i_1>i_2>\cdots>i_r\geq 1}\prod_{j=1}^r x_{i_j}^{2^{n+1-j+1-i_j}}.
\end{split}
\uqe
On the other hand,
\equ\label{proof_square_symetric_2}
\begin{split}
x_{n+1}S_{n,r-1}(x_1,\ldots,x_n)=x_{n+1}\sum_{n\geq i_1>i_2>\cdots>i_{r-1}\geq 1}\prod_{j=1}^{r-1} x_{i_j}^{2^{n-j+1-i_j}}\cr
\sum_{n\geq i_1>i_2>\cdots>i_{r-1}\geq 1}x_{n+1}\prod_{j=1}^{r-1}
x_{i_j}^{2^{n-j+1-i_j}}.
\end{split}
\uqe
Now, making $i_1=n+1$ and  $i_{j+1}=i_{j}$ (moving the variable $ j $ to $ j + 1$), we obtain that (\ref{proof_square_symetric_2}) becomes
\equ\label{proof_square_symetric_3}
\sum_{n+1= i_1>i_2>\cdots>i_r\geq 1}\prod_{j=1}^r x_{i_j}^{2^{n+1-j+1-i_{j}}}.
\uqe
Notice that the variable $x_{i_j}$ with exponent  $n-j+1-i_{j}$ in (\ref{proof_square_symetric_2}) concide with the variable $x_{i_{j+1}}$ with exponent  $n+1-j+1-i_{j+1}$ in (\ref{proof_square_symetric_3}).

Now, clearly the sum of (\ref{proof_square_symetric_1}) and (\ref{proof_square_symetric_3}) proves the lemma.
\frp

\pro\label{tees_formula}
For$$e_k(w)=\sum_{i=0}^{k-2}T_{k,i}[1]_w^{2^i},$$ is satisfied that
\begin{align*}
T_{k,i}=S_{k-2,k-2-i}(D_{2},D_{3},\ldots,D_{k-1}),
\end{align*}
where $D_i=e_i(t_i)$.
\orp

\prf
Using the identity $$e_{k+1}(w)=e^2_k(w)+D_ke_k(w),\mbox{ for $k\geq 2$},$$ we
obtain the following recursive equations
\begin{align*}
T_{k+1,0} &= D_kT_{k,0} \nonumber\\
T_{k+1,i} &=T^2_{k,i-1}+D_kT_{k,i}\\
T_{k+1,k-1} &=1\nonumber
\end{align*}
%
Then, from induction suppose that the proposition is valid for $T_{k,i}$, using the recursive form we get

\begin{align*}
T_{k+1,i}&= T^2_{k,i-1}+D_kT_{k,i}\cr
&=S^2_{k-2,k-2-(i-1)}(D_2,D_3,\ldots,D_{k-1})+D_kS_{k-2,k-2-i}(D_2,D_3,\ldots,D_{k-1})\cr
&=S^2_{k-2,k-1-i}(D_2,D_3,\ldots,D_{k-1})+D_kS_{k-2,k-1-i-1}(D_2,D_3,\ldots,D_{k-1})\cr
&=S_{k-1,k-1-i}(D_2,D_3,\ldots,D_{k}).
\end{align*}

The last equality follows from lemma (\ref{square_symetric}).
Now, the result follows from verifying that the coefficients  $T_{k,i}$ coincide with $S_{k-2,k-2-i}(D_2,D_3,\ldots,D_{k-1})$ for some first small values of $k$. 
\frp

For simplicity, set  $S_{k-2,k-2-i}:=S_{k-2,k-2-i}(D_{2},D_{3},\ldots,D_{k-1})$.

\crl\label{bees_formula}
The coefficients of the polynomial  $$e_k(w)=\sum_{i=0}^{k-1} B_{k,i} w^{2^i}$$
are given by the formulas
\begin{eqnarray*}
B_{k,k-1} & =&T_{k,k-2}=S_{k-2,0}=1\\
B_{k,i} &=&T_{k,i}+T_{k,i-1}=S_{k-2,k-2-i}+S_{k-2,k-1-i},\mbox{ for $1\leq i\leq k-2$}\\
B_{k,0} &=&T_{k,0}=S_{k-2,k-2}=D_{k-1}D_{k-2}\cdots D_2.
\end{eqnarray*}
\lrc
\prf
Note that
\begin{align*}
\begin{split}
e_k(w)&=\sum_{i=0}^{k-2}T_{k,i}[1]_w^{2^i}\cr
&=\sum_{i=0}^{k-2}T_{k,i}\left(w+w^2\right)^{2^i}\cr
&=T_{k,k-2}w^{2^{k-1}}+\sum_{i=1}^{k-2}\left(T_{k,i}+T_{k,i-1}\right)w^{2^i}+T_{k,0}w\mbox{ .}
\end{split}
\end{align*}
\frp

\section{Relationship among the values $d_k$, $\ell_k$ y $D_k$.}

Basically, these relationships are corollary of theorem (\ref{division_theorem}) and the explicit expression of the coefficients $B_{k,i}$ developed in the previous section.

If we evaluate the polynomial equality
\equ\label{polynomial_equality}
p_k(w)=\frac{e^2_k(w)}{d_k}+Ce_k(w)
\uqe
in $w=t_k$, we get that
$$1=\frac{D^2_k}{d_k}+C D_k.$$
Solving for $C$, we obtain
\equ\label{first_C_expresion}
C=\frac{1}{D_k}+\frac{D_k}{d_k}=\frac{d_k+D_k^2}{D_kd_k}.
\uqe
Now, using the definition of $C$ in  (\ref{division_theorem}),
we also have that
\begin{align*}
C=\frac{1}{d_{k-1}}+\frac{1+D^2_{k-1}+D_{k-2}^4+\cdots+D_2^{2^{k-2}}}{d_k},
\end{align*} since
$$B_{k,k-2}^2=(1+D_{k-1}+D_{k-2}^2+\cdots+D_2^{2^{k-3}})^2,$$
from corollary \ref{bees_formula} and part 2) of lemma
\ref{square_symetric}.

Multiplying by $D_k d_k$, we obtain
$$
CD_kd_k=\frac{D_kd_k}{d_{k-1}}+D_k\left(1+D^2_{k-1}+D_{k-2}^4+\cdots+D_2^{2^{k-2}}\right)$$
and using (\ref{first_C_expresion}), we have
$$
d_k+D_k^2=\frac{D_kd_k}{d_{k-1}}+D_k\left(1+D^2_{k-1}+D_{k-2}^4+\cdots+D_2^{2^{
k-2}}\right).
$$
From where,
$$d_k\left(1+\frac{D_k}{d_{k-1}}\right)=D_k(1+D_k+D^2_{k-1}+\cdots+D_2^{2^{k-2}}),$$
so eventually we get
\equ\label{recursive_formula_for_D_k}
\begin{split}
d_k&=\frac{D_kd_{k-1}}{d_{k-1}+D_k}\left(1+D_k+D^2_{k-1}+\cdots+D_2^{2^{k-2}}\right)\cr
&=\frac{D_kd_{k-1}}{d_{k-1}+D_k}\cdot B_{k+1,k-1}.
\end{split}
\uqe

Now, using the recursive formula (\ref{recursive_Dees}) is easy to see that
$$d_2=[1]_x$$
and also
$$D_2=e_2(t_2)=[1]_{t_2}=[1]_x,$$
equation (\ref{recursive_formula_for_D_k}) gives a recursive procedure to calculate $d_k$, in terms of values $D_i$'s  with $2\leq i \leq k$.

Now, equating the coefficients of the linear terms of the polynomials in (\ref{polynomial_equality}), we obtain that
$$\frac{1}{\ell_k}=CD_{k-1}D_{k-2}\cdots D_2$$
and using (\ref{first_C_expresion}), we conclude that
\begin{align*}
\ell_k=\frac{D_kd_k}{(d_k+D_k^2)(D_{k-1}D_{k-2}\cdots D_2)}.
\end{align*}

We summarize the above discussion in the main result of the article.

\begin{theorem}
 Recursive formulas to compute $\ell_k$ y $d_k$ values in terms of $D_k$'s.
\begin{enumerate}[1)]
\item $d_2=D_2$.
\item $d_k=\frac{D_kd_{k-1}}{d_{k-1}+D_k}\left(1+D_k+D^2_{k-1}+\cdots+D_2^{2^{k-2}}\right)$.
\item $\ell_k=\frac{D_kd_k}{(d_k+D_k^2)(D_{k-1}D_{k-2}\cdots D_2)}$.
\end{enumerate}
\end{theorem}

\subsubsection*{Acknowledgements}
Francisco Portillo thanks Conacyt for financial support for a year sabbatical stay at {\it Universidad Aut\'onoma de Yucat\'an} under the project grant \#261761.



\end{document}